\documentclass[a4]{article}

\usepackage{graphicx}
\usepackage{amsthm}

         \usepackage{verbatim}
         \usepackage[english]{babel}
         \usepackage{amsfonts,latexsym,amsmath,amssymb}

\newtheorem{dfn}{Definition}
\newtheorem{thm}{Theorem}
\newtheorem{pro}[thm]{Proposition}

\newtheorem{cor}[thm]{Corrolary}
\newcommand{\rr}{\mathbf r}
\newcommand{\eps}{\varepsilon}
\newcommand{\qq}{\quad}
\newcommand{\Ra}{\;\Rightarrow\;}
\newcommand{\lb}{[}
\newcommand{\rb}{]}

\newcommand{\D}{\displaystyle}

\newcommand{\ic}{{\,\mathcal I}}
\newcommand{\jc}{{\,\mathcal J}}
\newcommand{\lc}{{\,\mathcal L}}

\newcommand{\NN}{\mathbb N}
\newcommand{\RR}{\mathbb R}

\newcommand{\TT}{\mathbb T}

\newcommand{\dve}{{\mathbf 2}}
\newcommand{\oR}{\overline{R}}
\newcommand{\uR}{\underline{R}}
\newcommand{\oRp}{\overline{R}_\xi}
\newcommand{\uRp}{\underline{R}_\xi}

\newcommand{\be}{\begin{eqnarray*}}
\newcommand{\ee}{\end{eqnarray*}}
\newcommand{\bP}{\begin{proof}}
\newcommand{\eP}{\end{proof}}
\newcommand{\ed}{\end{document}}

\newcommand{\lf}{\lfloor}
\newcommand{\rf}{\rfloor}

\newcommand{\nas}[1]{\langle#1\rangle}
\newcommand{\mapp}{\longrightarrow}

\begin{document}
\begin{center}
  \vspace*{5\unitlength}
  {\LARGE {\bf Local return rates in Sturmian subshifts}}\\
  \bigskip
  Michal Kupsa\\
  \bigskip
  {\it Faculty of Mathematics and Physics, Charles University in Prague\\
    and\\
    Centre de Physique Th\'eorique, CNRS Luminy, Marseille}\\
  \bigskip
  \begin{minipage}[h]{0.8\linewidth}
    \small The local return rates have been introduced by Hirata,
    Saussol and Vaienti \cite{HSV99} as a tool for the study of the
    asymptotic distribution of the return times to cylinders. We give
    formulas for these rates in Sturmian subshifts.
  \end{minipage}
\end{center}
\bigskip

\section{Introduction}
The lower and upper local return rates have been introduced by Hirata,
Saussol and Vaienti in \cite{HSV99} as a tool for the study of the
asympotic distribution of the return times to cylinders in a class of
non-uniformly hyperbolic dynamical systems. They are functions
$\uRp,\oRp:X\mapp \lb 0,\infty]$ defined for an arbitrary topological
dynamical system $(X,F)$ and a finite partition $\xi$ of $X$. For a
subshift $\Sigma\subseteq A^\NN$ and the canonical partition $\{\lb
a]\mid a\in A\}$ we can reformulate the definition as \be
\uR(x) &=& \liminf_{n\to\infty}\frac{\tau\left(\lb x(n)]\right)}{n}\\
\oR(x) &=& \limsup_{n\to\infty}\frac{\tau\left(\lb x(n)]\right)}{n}.
\ee Here $x(n)=x_0x_1\ldots x_{n-1}$ is a prefix of $x\in\Sigma$ of
length $n$, $\lb x(n)]$ is its cylinder and $\tau(\lb x(n)])$ is the
Poincar\'e return time of $\lb x(n)]$.

For an arbitrary dynamical system $(X,F)$ the functions $\uRp,\oRp$
are subinvariant, i.e., $\uRp\circ F\leq\uRp$ and $\oRp\circ
F\leq\oRp$.  Moreover, if $\mu$ is an $F$-invariant Borel probability
measure and $\xi$ is a measurable partition of $X$, then $\uR$ and
$\oR$ are invariant allmost everywhere. In particular, if $(X,F,\mu)$
is ergodic, then by the Birkhoff ergodic theorem there exist constants
$\rr_0,\rr_1\in\lb 0,\infty]$ such that for almost all $x\in X$,
$\uRp(x)=\rr_0$ and $\oRp(x)=\rr_1$.

The ergodic case has been treated in several more papers. Saussol et
al \cite{STV02}( see also \cite{ACS03}) show that if the entropy of
$\mu$ is positive, then $\rr_0 \geq 1$. Cassaigne et al \cite{CHV}
show that this inequality is not satisfied for systems with zero
entropy. In particular for the Fibonacci shift obtained from the
golden angle rotation, the lower local rate assumes the value
$\rr_0=\frac{3-\sqrt{5}}{2}<1$. Afraimovich et al \cite{ACS03} show
that $\rr_0 = 0$ for some rotations of the circle whose parameter has
unbounded continued fraction expansion. It follows that the same
result holds for the corresponding Sturmian subshift. K\r{urka}
\cite{Ku03} treats the case of substitutive subshifts and obtains a
formula for $\rr_0$ and $\rr_1$. In this case both $\rr_0$ and $\rr_1$
are positive and finite.

In this paper we will discuss completely the situation in the Sturmian
shifts.  One can easy check that the result of Afraimovich et al
considered for corresponding Sturmian shifts and the result of
Cassaigne et al for Fibonacci shift follows immediately.  We give
formulas for $\rr_0$ and $\rr_1$ in terms of the convergents $q_k$
obtained from the continued fraction expansion of the parameter
$\alpha = \lb 0,a_1,a_2,\ldots]$. If $a_k$ are bounded, then $\rr_0$
and $\rr_1$ are positive and finite. If $a_k$ are unbounded, then
$\rr_0=0$ and $\rr_1 = \infty$. This result, that $\rr_0=0$ iff
$\rr_1=\infty$ iff the continued fraction expansion is unbounded, has
been obtained by a different technique by Chazottes and Durand in
\cite{CD05}

\section{Sturmian shifts}

A dynamical system is a pair $(X,F)$, where $X$ is a compact metric
space and $F$ is a continuous function from $X$ to $X$. The Poincar\'e
return time of a subset $M\subseteq X$ is
$$\tau(M)=\min\{k>0\mid F^k(M)\cap M\neq\emptyset\}.$$
Let $A$ be a finite alphabet, and $A^\NN$ the space of all infinite
sequences of letters from $A$ with the product topology. The set $A^*$
consists of all words (finite sequences) over $A$. For a word
$u=u_0u_2\ldots u_{n-1}\in A^*$, denote by $|u|=n$ its length.  The
set $A^n$ consists of all words of length $n$. The shift map
$\sigma:A^\NN\longrightarrow A^\NN$ is defined by
$\sigma_i(x)=x_{i+1}$.

A shift is any subsystem $(\Sigma,\sigma)$ of $(A^{\NN},\sigma)$,
where $\Sigma\subseteq A^\NN$ is nonempty, closed and
$\sigma$-invariant.  For a shift $\Sigma$ and for a word
$u=u_0u_1\ldots u_{n-1}\in A^*$ we denote by $\lb
u]=\{x\in\Sigma\mid\forall i<n: x_i=u_i\}$ the cylinder of $u$. The
language of a shift is the set of words which have nonempty cylinders,
i.e., $\lc(\Sigma)=\{u\in A^*\mid \lb u]\neq\emptyset\}$. The set
$\lc^n(\Sigma)$ consists of all words of the language of length $n$.
If we denote by $x(n)=x_0x_1\ldots x_{n-1}$ the prefix of $x\in\Sigma$
of length $n$, then $\lc^n(\Sigma)=\{x(n)\mid x\in\Sigma\}$.

A Sturmian shift is a coding of an irrational rotation of the unit
circle (Hedlund and Morse \cite{HM40}). This is a dynamical system
$(\TT,F_{\alpha})$, where $\TT=\lb 0,1\lb $ is the circle with the
metric $d(x,y)=\min\{|x-y|,1-|x-y|\}$ and $F_{\alpha}(x)=x+\alpha\mod\
1$, where $\alpha\in\RR$. We consider only irrational angles from the
open interval $\alpha \in ]0,1\lb $.

There is the canonical partition $\ic = \{I_0,I_1\}$ of $\TT$, where
$I_0=\lb 0,1-\alpha\lb $ and $I_1=\lb 1-\alpha,1\lb $. For
$u\in\dve^*$, set
$$I_u=\bigcap^{|u|-1}_{k=0} F_{\alpha}^{-k}(I_{u_k}).$$
Any $I_u$ is either a semiopen interval or the empty set. The
associated Sturmian shift $(\Sigma_{\alpha},\sigma)$ is defined by its
language $\lc(\Sigma_{\alpha})=\{u\in\dve^*\mid I_u\neq\emptyset\}$.
In other words,
$$\Sigma_{\alpha}=\{x\in\dve^\NN\mid
\forall n\in\NN,I_{x(n)}\neq\emptyset\}.$$ If $\alpha\in]0,1\lb $ is
irrational, both the rotation $(\TT,F_{\alpha})$ and the Sturmian
shift $(\Sigma_{\alpha},\sigma)$, are minimal and uniquely
ergodic. Moreover, if $u\in\lc(\Sigma_{\alpha})$, then
$$\mu(\lb u])=|I_u|,\qq \tau(\lb u])=\tau(I_u),$$
where $|I_u|$ is the length of the interval $I_u$.  It follows that
the local return rates can be computed from the return times of
intervals.

\be
\uR(x)&=&\liminf_{n\to\infty}\frac{\tau(I_{x(n)})}{n}\\
\oR(x)&=&\limsup_{n\to\infty}\frac{\tau(I_{x(n)})}{n}.  \ee

The description of the intervals $I_u$ is obtained from the continued
fraction expansion of $\alpha$. There exists a unique sequence
$\{a_k\}^\infty_{k=1}$ of positive integers such that
$$\alpha=\lb 0,a_1,a_2,\ldots]=0 + \frac{1}{\D a_1+\frac{1}{\D a_2+\ldots}}.$$
The convergents of $\alpha$ are the sequences $\{p_k\}^\infty_{k=-1}$,
$\{q_k\}^\infty_{k=-1}$ defined by $p_{-1}=1$, $q_{-1}=0$, $p_0=0$,
$q_0=1$ and
$$q_{k+1}=a_{k+1}q_k+q_{k-1},\qq
p_{k+1}=a_{k+1}p_k+p_{k-1}.$$ By the Klein theorem (see Hardy and
Wright \cite{HW80}), the closest returns of the iterates
$F_{\alpha}^n(0)$ to zero happen at times $q_k$. We have
$d(0,F^{q_k}(0)) = \eta_k = (-1)^{k}(q_k\alpha-p_k)$ and for
$q_k<n<q_{k+1}$, $d(0,F^n(0))>\eta_k$.  In particular $\eta_{-1}=1$,
$\eta_0=\alpha$ and
$$ \eta_{k+1}=a_{k+1}\eta_k - \eta_{k-1}.$$
The sequence $\{\eta_k\}^\infty_{k=-1}$ is positive, decreasing and
converges to zero. It follows that if $I=\lb a,b\lb $ is a semiopen
interval, then
$$ \eta_{k+1}<|I|\leq\eta_k \Ra \tau(I)=q_{k+1}.$$
The return times of intervals from $\ic^n$ are therefore
convergents $q_k$. We determine times when the return times jump
from some $q_k$ to a higher $q_{k+1}$ (or $q_{k+2}$) and obtain a
formula for the local return rates.

\section{Jumps of the return time}

\begin{pro}
  For $x\in\Sigma_{\alpha}$, $k \geq -1$, define the $k$-th jump of
  the return time as
$$r_k(x)=\min\{n\in\NN\mid\tau(I_{x(n)})\geq q_{k+1}\}.$$
Then $r_{-1}(x)=0$ and the following equalities hold for $x\in
\Sigma_{\alpha}$.  \be \uR(x)&=&\liminf_{k\to\infty}\frac{q_k}{r_k(x)}
=
1/\limsup_{k\to\infty}\frac{r_k(x)}{q_k} \\
\oR(x)&=&\limsup_{k\to\infty}\frac{q_{k+1}}{r_k(x)} =
1/\liminf_{k\to\infty}\frac{r_k(x)}{q_{k+1}} \ee \end{pro}

\bP
For $x\in\Sigma_{\alpha}$, denote $S=\{k\in\NN\mid
r_{k-1}(x)<r_k(x)\}$. The set is infinite and we can order it into
increasing sequence $\{k_i\}^\infty_{i=0}$.  If $r_{k_i}(x)\leq
n<r_{k_{i+1}}(x)$, then $\tau(I_{x(n)})=q_{k_{i+1}}$ and if
$k_i<k<k_{i+1}$, then
$\frac{q_k}{r_k(x)}\geq\frac{q_{k_i}}{r_{k_i}(x)}$,
$\frac{q_k}{r_{k-1}(x)}\leq\frac{q_{k_{i+1}}}{r_{k_i}(x)}$. Thus \be
\uR(x)&=&\liminf_{n\to\infty}\frac{\tau(I_{x(n)})}{n}=
\liminf_{i\to\infty}\left(\min_{r_{k_i}(x)\leq n<r_{k_{i+1}}(x)}
  \frac{\tau(I_{x(n)})}{n}\right)=\\
&=&\liminf_{i\to\infty}\frac{q_{k_{i+1}}}{r_{k_{i+1}}(x)-1}=
\liminf_{i\to\infty}\frac{q_{k_i}}{r_{k_i}(x)}=
\liminf_{k\to\infty}\frac{q_k}{r_k(x)}.\\
\oR(x)&=&\limsup_{i\to\infty}\left(\max_{r_{k_i}(x)\leq
    n<r_{k_{i+1}}(x)} \frac{\tau(I_{x(n)})}{n}\right)=
\limsup_{k\to\infty}\frac{q_{k+1}}{r_k(x)}.  \ee \eP

To compute the jumps of the return time, we construct another
symbolic description of Sturmian shifts.  The partition
$\ic^n=\{I_u\mid u\in\lc^n(\Sigma_{\alpha})\}$ consists of semiopen
intervals on the unit circle divided by cut points
$$ \mbox{Cut}(n)=\{\nas{i}\mid i=0,1,\ldots,n\},$$
where $\nas{i}= F_{\alpha}^{-i}(0) = (-i\alpha)\mod 1$.  The structure
of $\ic^n$ is described by the Three length theorem (S\'os
\cite{So58}) which says that $\ic^n$ contains intervals of at most
three lengths. For some $n$, however $\ic^n$ contains only intervals
of two lengths. This happens in particular at times $n=q_k-1$, when
the intervals of $\ic^n$ have lengths $\eta_{k-1}$ and
$\eta_{k-1}+\eta_k$. To describe the partitions $\ic^{q_k-1}$ we
consider a new symbolic space $X_{\alpha}$ which consists of paths in
the infinite graph in Figure \ref{graph}.  It looks like Bratelli
diagram( \cite{DDM00}), but the dynamics on $X_{\alpha}$ is far more
complicated. The main reason for introducing the space $X_{\alpha}$ is
to obtain a simple formula for $r_k(x)$ in Proposition \ref{rk}.

\setlength{\unitlength}{0.8mm}
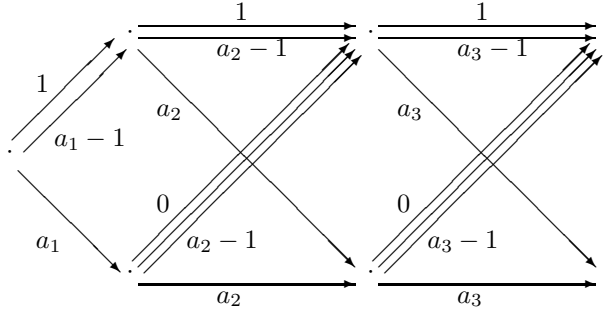
\begin{figure}
  \begin{picture}(100,50) \put(0,25){$\cdot$} \put(20,5){$\cdot$}
    \put(20,45){$\cdot$} \put(60,5){$\cdot$} \put(60,45){$\cdot$}
    \put(100,5){$\cdot$} \put(100,45){$\cdot$}
    \put(1,28){\vector(1,1){17}} \put(3,26){\vector(1,1){17}}
    \put(2,23){\vector(1,-1){17}} \put(22,47){\vector(1,0){36}}
    \put(22,45){\vector(1,0){36}} \put(22,43){\vector(1,-1){36}}
    \put(22,4){\vector(1,0){36}} \put(21,8){\vector(1,1){36}}
    \put(22,7){\vector(1,1){36}} \put(23,6){\vector(1,1){36}}
    \put(62,47){\vector(1,0){36}} \put(62,45){\vector(1,0){36}}
    \put(62,43){\vector(1,-1){36}} \put(62,4){\vector(1,0){36}}
    \put(61,8){\vector(1,1){36}} \put(62,7){\vector(1,1){36}}
    \put(63,6){\vector(1,1){36}} \put(5,36){$1$} \put(8,27){$a_1-1$}
    \put(5,10){$a_1$} \put(38,48){$1$} \put(35,42){$a_2-1$}
    \put(25,32){$a_2$} \put(25,16){$0$} \put(30,10){$a_2-1$}
    \put(35,1){$a_2$} \put(78,48){$1$} \put(75,42){$a_3-1$}
    \put(65,32){$a_3$} \put(65,16){$0$} \put(70,10){$a_3-1$}
    \put(75,1){$a_3$}
  \end{picture}
  \caption{The symbolic space $X_{\alpha}$ \label{graph}}
\end{figure}

\begin{dfn}
  For an irrational $\alpha=\lb 0,a_1,a_2,\ldots]$ set \be X_{\alpha}
  &=& \left\{x\in \prod_{k=1}^{\infty} \{0,1,\ldots,a_k\}\mid
    x_1\neq 0,\qq (x_{k+1}=0 \Ra x_k=a_k) \right\} \\
  \lc^n(X_{\alpha}) &=& \left\{u\in \prod_{k=1}^{n}
    \{0,1,\ldots,a_k\}\mid
    u_1\neq 0,\qq (u_{k+1}=0 \Ra u_k=a_k) \right\} \\
  \lc(X_{\alpha}) &=& \bigcup_{n\geq 1} \lc^n(X_{\alpha}).  \ee We
  construct a system of intervals $\{J_u\mid u \in \lc(X_{\alpha})\}$.
  If $1\leq u_1 \leq a_1$ set
$$  J_{u_1} = \left\{ \begin{array}{l@{\quad\mbox{if}\quad}l}
    \lb  \nas{u_1},\nas{u_1-1} \lb  & u_1<a_1 \\
    \lb  \nas{0}, \nas{a_1-1} \lb   & u_1=a_1 \end{array} \right. $$
If $u\in \lc^k(X_{\alpha})$, $k>1$, $J_{u(k-1)} = (-1)^{k-2}\lb \nas{a},\nas{b}\lb $,
and if $u_{k-1}<a_{k-1}$  set
$$ J_u = \left\{\begin{array}{l@{\quad\mbox{if}\quad}l}
    (-1)^{k-1} \lb \nas{u_kq_{k-1}+a},\nas{(u_k-1)q_{k-1}+a}\lb  &
    1\leq u_k\leq a_k-1  \\
    (-1)^{k-1}\lb \nas{b},\nas{(a_k-1)q_{k-1}+a}\lb    &
    u_k=a_k  \end{array} \right. $$
If $u_{k-1}=a_{k-1}$ set
$$ J_u = \left\{\begin{array}{l@{\quad\mbox{if}\quad}l}
    (-1)^{k-1}\lb \nas{(u_k+1)q_{k-1}+a},\nas{u_kq_{k-1}+a}\lb  &
    0\leq u_k\leq a_k-1 \\
    (-1)^{k-1}[ \nas{b},\nas{u_kq_{k-1}+a}[  &
    u_k=a_k \end{array}\right. $$
\end{dfn}

Here $(-1)\lb b,a\lb =\lb a,b\lb $, where $0 \leq a<b<1$, is a
semiopen interval of the circle. We identify also $\lb a,0\lb =\lb
a,1\lb =(-1)\lb 0,a\lb $. If $x \in X_{\alpha}$, we denote by $x(n) =
x_1\ldots x_n$ the prefix of $x$ of length $n$.  In Figure \ref{part}
we can see the partitions of the circle for $\alpha=\lb
0,2,3,\ldots]$.

\setlength{\unitlength}{1.2mm}
\begin{figure}
  \begin{picture}(100,30) \put(-1,4.5){\line(1,0){102}}
    \put(0,4){\line(0,1){1}} \put(13,4){\line(0,1){1}}
    \put(26,4){\line(0,1){1}} \put(39,4){\line(0,1){1}}
    \put(56.5,4){\line(0,1){1}} \put(69.5,4){\line(0,1){1}}
    \put(82.5,4){\line(0,1){1}} \put(100,4){\line(0,1){1}}
    \put(2,2){$I_{001010}$} \put(15,2){$I_{010010}$}
    \put(28,2){$I_{010100}$} \put(43,2){$I_{010101}$}
    \put(59,2){$I_{100101}$} \put(72,2){$I_{101001}$}
    \put(87,2){$I_{101010}$} \put(-2,0){$\nas{0}$}
    \put(11,0){$\nas{2}$} \put(24,0){$\nas{4}$} \put(37,0){$\nas{6}$}
    \put(54.5,0){$\nas{1}$} \put(67.5,0){$\nas{3}$}
    \put(80.5,0){$\nas{5}$} \put(5,6){$J_{20}$} \put(18,6){$J_{21}$}
    \put(31,6){$J_{22}$} \put(47,6){$J_{23}$} \put(62,6){$J_{11}$}
    \put(75,6){$J_{12}$} \put(88,6){$J_{13}$}
    \put(-1,14.5){\line(1,0){102}} \put(0,14){\line(0,1){1}}
    \put(56.5,14){\line(0,1){1}} \put(100,14){\line(0,1){1}}
    \put(30,12){$I_0$} \put(80,12){$I_1$} \put(30,16){$J_2$}
    \put(80,16){$J_1$} \put(-2,11){$\nas{0}$} \put(54.5,11){$\nas{1}$}
  \end{picture}
  \caption{Partitions of the circle \label{part}}
\end{figure}
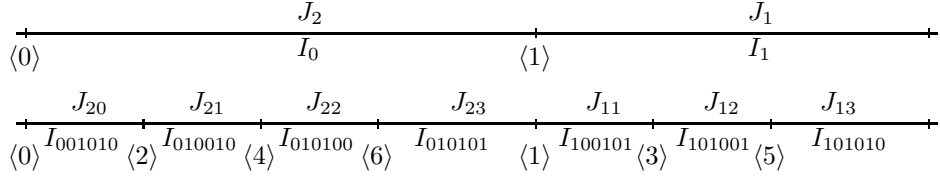

\begin{pro} \label{intervals} If $u \in \lc^k(X_{\alpha})$, $k\geq 1$,
  then
$$ |J_u| = \left\{ \begin{array}{l@{\quad\mbox{if}\quad}l}
    \eta_{k-1} & u_{k}<a_k \\ \eta_{k-1}+\eta_k & u_k=a_k
  \end{array} \right.$$
and
\be
\ic^{q_k-1} &=& \{J_u\mid u \in \lc^k(X_{\alpha})\} \\ &=&
\{(-1)^{k-1}\lb \nas{i+q_{k-1}},\nas{i}\lb  \mid
i=0,1,\ldots,q_k-q_{k-1}-1\} \cup\\
&&\{(-1)^{k-1}\lb \nas{i},\nas{i+q_k-q_{k-1}}\lb  \mid
i=0,1,\ldots,q_{k-1}-1\}.
\ee
Moreover, $\jc_u=\{J_v\mid v\in \lc^{k+1}(X_{\alpha})| v(k)=u\}$
is a partition of $J_u$ and
$$\ic^{q_{k+1}-1}=\bigcup_{u\in\lc^k(X_{\alpha})}\jc_u.$$
\end{pro}

\bP If $u\in\lc^k(X_\alpha)$, $k\geq 1$, $u_k<a_k$, then $J_u$ is an
image of $(-1)^{k-1}\lb \nas{q_{k-1}},0\lb $ in a rotation. By the
Klein theorem, $|I_u|=\eta_{k-1}$.

We have $1>\nas{1}>\nas{2}> \cdots \nas{a_1-1}>0$, so $\{J_{u_1}\mid
1\leq u_1\leq a_1\} = \ic^{q_1-1}$ and $|J_{a_1}| = 1-(a_1-1)\eta_0 =
\eta_0+\eta_1$.  Assume that the first part of the proposition holds
for $k\geq 1$.  Let $u \in \lc^k(X_{\alpha})$. Intervals from
$M=\{J_{uj}\mid j<a_{k+1}\}$ coincide and we have proved that its
length is $\eta_k$.  Denote $J=\bigcup M$. If $u_k<a_k$ then
$|J|=(a_{k+1}-1)\eta_k$ and if $u_k=a_k$ then $|J|=a_{k+1}\eta_k$. In
both cases, $|J|<|J_u|$, $J_{ua_{k+1}}=J_u-J$ and
$|J_{ua_{k+1}}|=\eta_{k-1}-(a_{k+1}-1)\eta_k=\eta_k+\eta_{k+1}$.  Thus
$\jc_u$ is a partition of $J_u$. Because $\{J_u\mid
u\in\lc^k(X_\alpha)\}$ is a partition of $\TT$, then also
$$\jc=\{J_v\mid v\in\lc^{k+1}(X_\alpha)\}=
\bigcup_{u\in\lc^k(X_\alpha)}\jc_u$$ is. It is not difficult to prove
that the endpoints of intervals from $\jc$ belong to $\{\nas{i}\mid
0\leq i\leq q_{k+1}-1\}$. The partitions $\jc$ and $\ic^{q_{k+1}-1}$
contain intervals of two lengths $\eta_k$ and $\eta_k+\eta_{k+1}$,
hence $\jc=\ic^{q_{k+1}-1}$. For the partition \be \jc'&=& \{(-1)^k\lb
\nas{i+q_k},\nas{i}\lb \mid
i=0,1,\ldots,q_{k+1}-q_k-1\} \cup\\
&&\{(-1)^k\lb \nas{i},\nas{i+q_{k+1}-q_k}\lb \mid i=0,1,\ldots,q_k-1\}
\ee we prove the equality $\ic^{q_{k+1}-1}=\jc'$ similarly.  \eP

For each $k\geq 1$ we have thus a one-to-one map $\gamma_k:
\lc^{q_k-1}(\Sigma_{\alpha}) \to \lc^k(X_{\alpha})$ given by
$J_{\gamma_k(u)} = I_{u}$. For the corresponding symbolic spaces we
get a homeomorphism $\gamma: \Sigma_{\alpha} \to X_{\alpha}$ given by
$\gamma(x)(k) = \gamma_k(x(q_k-1))$. The local return rates, as well
as the functions of the return jumps are carried over to the space
$X_{\alpha}$. By the abuse of notation we keep for them the same
symbols $\underline{R},\overline{R}: X_{\alpha} \to \lb 0,\infty]$,
$r_k:X_{\alpha} \to \NN$. We now obtain a recursive formula for $r_k$.

\begin{pro} \label{rk} For $x \in X_{\alpha}$ we have $r_{-1}(x) = 0$
  and
$$ r_k(x) = x_{k+1}q_k + r_{k-1}(x) = \sum_{j=0}^k x_{j+1}q_j. $$
\end{pro}

\bP Assume $y \in \Sigma_{\alpha}$, $x=\gamma(y) \in X_{\alpha}$ and
$k\geq 0$. We show first that if $J_{x(k)} = (-1)^{k-1}\lb
\nas{a},\nas{b}\lb $, then $r_{k-1}(x)=b+q_{k-1}$. If $x_{k}<a_{k}$,
then $J_{x(k)} = I_{y(a)}$.  Since $I_{y(a-1)} \neq I_{y(a)}$,
$|I_{y(a-1)}|>\eta_{k-1}$ and $r_{k-1}(x) = a = b+q_{k-1}$.
Let $x_{k}=a_{k}$. Since the form of partition $\{J_u\mid
u\in\lc^{k+1}(X_\alpha),u(k)=x(k)\}$ of $J_{x(k)}$ we get $J_{x(k)}=
I_1\cup I_2$ where
$$I_1=(-1)^{k-1}\lb \nas{a},\nas{a+q_{k}}\lb ,\qq
I_2=(-1)^{k-1}\lb \nas{a+q_{k}},\nas{a+q_{k}-q_{k-1}}\lb ,$$
$I_1,I_2\in\ic^{a+q_{k}}$, $|I_1|=\eta_{k}$, $|I_2|=\eta_{k-1}$ and
$I_{y(a+q_k-1)}=J_{x(k)}$ and either $I_{y(a+q_k)}=I_1$ or
$I_{y(a+q_k)}=I_2$. Hence $|I_{y(a+q_k)}|\leq\eta_{k-1}$,
$|I_{y(a+q_k-1)}|>\eta_{k-1}$ and $r_{k-1}(x)=a+q_k=b+q_{k-1}$.

Assume now that $J_{x(k+1)} = (-1)^k\lb \nas{c},\nas{d}\lb $, so
$r_k(x)=d+q_k$. Put j=1 if $x_k=a_k$, j=0 otherwise. It follows
$d=a+(x_{k+1})-1)q_k+jq_k$ and $a=b+q_{k-1}-jq_k$. Thus \be
r_k(x)-r_{k-1}(x)&=&(d+q_k)-(b+q_{k-1})=(a+x_{j+1}q_k+jq_k)-
(a+jq_k)\\
&=& x_{k+1}q_k.  \ee \eP

\begin{pro} \label{boundrk} For every $x\in X_{\alpha}$ we have $q_k
  \leq r_k(x) \leq q_{k+1}+q_k-1$.
\end{pro}

\bP Clearly $q_{-1} = 0 = r_{-1}(x) = 0 = q_0 + q_{-1} -1$, $q_1 = 1
\leq r_1(x) \leq a_1 = q_1 +q_0 -1$. Assume that the statement holds
for all integers less than $k$. Then
$$ r_k(x) = x_{k+1}q_k +r_{k-1}(x) \leq a_{k+1}q_k + q_k +q_{k-1} -1
= q_k + q_{k+1} -1 $$ If $x_{k+1} \geq 1$, then $r_k(x) = x_{k+1}q_k
+r_{k-1}(x) \geq q_k$.  If $x_{k+1}=0$, then $x_k=a_k$ and
$$ r_k(x) = r_{k-1}(x) = a_kq_{k-1}+r_{k-2}(x)
\geq a_kq_{k-1} +q_{k-2} = q_k. $$ \eP

\begin{pro}\label{meze}
  Define the points $b,c,d\in X_{\alpha}$ by
$$ b = (a_1,a_2,a_3,\ldots),\qq
c = (1,a_2,0,a_4,0,a_6,\ldots),\qq d = (a_1,0,a_3,0,a_5,\ldots ) $$
Then
\begin{eqnarray*}
  \min \uR &=& \uR(b) \;=\; \liminf_{k\to\infty}\frac{q_k}{q_{k+1}+q_k-1}
  \;=\; \rr_0  \\
  \min \oR &=& \oR(b) \;=\; \limsup_{k\to\infty}\frac{q_{k+1}}{q_{k+1}+q_k-1}\\
  \max \oR&=&\max(\oR(c),\oR(d))=
  \limsup_{k\to\infty}\frac{q_{k+1}}{q_k} \;=\; \rr_1.
\end{eqnarray*}
\end{pro}
\bP It is easy to see that for $k\in\NN$, \be
r_k(b)&=&\sum_{j=0}^k a_{j+1}q_j=q_{k+1}+q_k-1\\
r_{2k-1}(c)=r_{2k}(c)&=&1+\sum_{j=1}^k a_{2j}q_{2j-1}=q_{2k}\\
r_{2k}(d)=r_{2k+1}(d)&=&\sum_{j=0}^k a_{2j+1}q_{2j}=q_{2k+1} \ee

By Proposition \ref{boundrk} we obtain the bounds for the limits in
the right hand sides. The following formulas complete the proof.
\begin{eqnarray*}
  \uR(b)&=&\liminf_{k\to\infty}\frac{q_k}{r_k(b)} = \rr_0\\
  \oR(b)&=&\limsup_{k\to\infty}\frac{q_{k+1}}{r_k(b)}=
  \limsup_{k\to\infty}\frac{q_{k+1}}{q_{k+1}+q_k-1}\\
  \max(\oR(c),\oR(d))&\geq&\max\left(
    \limsup_{k\to\infty}\frac{q_{2k+1}}{r_{2k}(c)},
    \limsup_{k\to\infty}\frac{q_{(2k-1)+1}}{r_{2k-1}(d)}\right)\\
  &\geq&\max\left(\limsup_{k\to\infty}\frac{q_{2k+1}}{q_{2k}},
    \limsup_{k\to\infty}\frac{q_{2k}}{q_{2k-1}}\right)\\
  &\geq&\limsup_{k\to\infty}\frac{q_{k+1}}{q_k} = \rr_1.
\end{eqnarray*}
\eP We have not been able to obtain a formula for $\min
\overline{R}$. Our results, however are sufficient to get formulas for
$\rr_0$ and $\rr_1$. Now, put some bounds for the values $\rr_0$ and
$\rr_1$.
\begin{pro}
  Let $\alpha=\lb 0,a_1,a_2,..]$ be irrational, $M=\limsup a_k$,
  $\gamma=\frac{\sqrt{5}+1}{2}$.  If the continued fraction expansion
  of $\alpha$ is unbounded, then $\rr_0=0$ and $\rr_1=\infty$.
  Otherwise, $\rr_1=\frac{1}{\rr_0}-1$ and
$$ \frac{1}{M+2} \leq \rr_0\leq\gamma^{-2} <\gamma\leq\rr_1\leq M+1$$
Moreover, $\rr_1=\gamma$( resp. $\rr_0=\gamma^{-2}$) if and only if
$M=1$.
\end{pro}
\bP Let $\alpha=\lb 0,a_1,a_2,..]$ be irrational, $M=limsup\ a_k$,
$\gamma=\frac{\sqrt{5}+1}{2}$. Denote $B_k=\frac{q_{k+1}}{q_k}$. Then
$a_k\leq B_k\leq a_{k+1}+1$ and $\rr_1=\limsup B_k$,
$$\rr_0=\liminf\frac{1}{B_k+1-\frac{1}{q_k}}=\liminf\frac{1}{B_k+1}$$
\indent If the continued fraction expansion of $\alpha$ is unbounded,
then also $\{B_k\}^\infty_{k=0}$ is. Hence $\rr_0=0$ and
$\rr_1=\infty$.

Let $M\in\NN$. Then $M\leq\rr_0\leq M+1$. If $M\geq 2$, then
$\gamma<M\leq\rr_0$. If $M=1$, then there exists $n_0\in\NN$, such
that for every $n>n_0$, $a_n=1$. Hence \be \limsup
B_k=1+\frac{1}{\liminf B_k}&,& \liminf B_k=1+\frac{1}{\limsup B_k} \ee
It implies that $\rr_0=1+\frac{1}{1+\frac{1}{\rr_0}}$. This equality
have just one positive solution $\rr_0=\gamma$. All properties of
$\rr_1$ is given by the equality $\rr_0=\frac{1}{\rr_1+1}$.  \eP

\section{The measure}

We are going to show that the constants $\rr_0$ and $\rr_1$ are
assumed by $\underline{R}$ and $\overline{R}$ almost everywhere.  The
unique invariant measure $\mu$ on $\Sigma_{\alpha}$ is carried over to
the space $X_{\alpha}$ using the length of associated intervals. If $u
\in \lc^k(X_{\alpha})$, then the measure of the cylinder of $u$ is
$\mu\{x\in X_{\alpha}| x(k)=u\} = |J_u|$. Define the projections $W_k:
X_{\alpha} \to \{0,1,\ldots,a_k\}$ by $W_k(x) = x_k$.  Then $W_k$ are
random variables and $(W_k)_{k\geq 1}$ is a nonstationary Markov
chain. Using Proposition \ref{intervals} we get the transition
probabilities.
$$ \begin{array}{rcl@{\quad\mbox{for}\quad}l}
  \mu \lb  W_1=j \rb &=& \eta_0=\alpha & 1\leq j<a_1 \medskip\\
  \mu \lb  W_1=j \rb &=& \eta_0+\eta_1 & j=a_1 \medskip\\
  \mu \lb  W_{k+1}=j|W_k<a_k \rb &=& \D\frac{\eta_k}{\eta_{k-1}} &
  1\leq j<a_{k+1} \medskip\\
  \mu \lb  W_{k+1}=j|W_k<a_k \rb &=& \D\frac{\eta_k+\eta_{k+1}}
  {\eta_{k-1}} &  j=a_{k+1} \medskip\\
  \mu \lb  W_{k+1}=j|W_k=a_k \rb &=& \D\frac{\eta_k}
  {\eta_{k-1}+\eta_k} &  0\leq j<a_{k+1} \medskip\\
  \mu \lb  W_{k+1}=j|W_k=a_k \rb &=& \D\frac{\eta_k+\eta_{k+1}}
  {\eta_{k-1}+\eta_k} &  j=a_{k+1} \end{array} $$

\begin{thm} \label{unbounded} If the continued fraction expansion is
  unbounded, then $\underline{R}(x) = 0$, $\overline{R}(x) = \infty$
  almost everywhere.
\end{thm}

\bP For every $x\in X_{\alpha}$ we have
$$ x_{k+1} \leq \frac{x_{k+1}q_k+r_{k-1}(x)}{q_k} = \frac{r_k(x)}{q_k}$$
Given $m\geq 1$ then $C_m = \{k\geq 1| a_k\geq m\}$ is an infinite
set.  Assume that $k+1 \in C_{2m+1}$. We have \be \mu\lb W_{k+1}\leq
m|W_k<a_k] &=& \frac{m\eta_k}{\eta_{k-1}} =
\frac{m\eta_k}{a_{k+1}\eta_k+\eta_{k+1}} \leq
\frac{m}{a_{k+1}} \leq \frac{1}{2} \\
\mu\lb W_{k+1}\leq m|W_k=a_k] &=&
\frac{(m+1)\eta_k}{\eta_{k-1}+\eta_k} \leq \frac{m+1}{a_{k+1}+1} \leq
\frac{1}{2} \ee It follows that $\mu\lb W_{k+1} \leq m|W_j=i] \leq
\frac{1}{2}$ for any $j\leq k$ and any $i\in \{0,1,\ldots,a_j\}$.
Given $k_0>0$ let $k_0<k_1 < \cdots k_n$ be a sequence of integers
from $C_{2m+1}$. Then \be \mu \lb W_{k_1}\leq m,\ldots,W_{k_n}\leq m]
&=& \mu\lb W_{k_1}\leq m]\cdot \mu\lb W_{k_2}\leq m|W_{k_1}\leq m]\cdots\\
&& \mu\lb W_{k_n}\leq m|W_{k_1}\leq m,\ldots, W_{k_{n-1}}\leq m]\\
&\leq& 2^{-n+1} \ee It follows
$$ \mu\left\{x\in X_{\alpha} \left| \frac{r_{k_i}(x)}{q_{k_{i}}}
    \leq m, 1\leq i\leq n \right \} \right. \leq \mu \{x\in X_{\alpha}
\mid x_{k_i} \leq m, 1\leq i \leq n\} \leq 2^{-n+1} $$ so $\mu\{x\in
X_{\alpha}| \underline{R}(x)>\frac{1}{m}\} = 0$ and $\underline{R}(x)
= 0$ almost everywhere.  We prove now the statement for
$\overline{R}$.  Given $\eps \in ]0,1\lb $, let $m$ be an integer with
$1-\eps+\frac{4}{m} = \delta<1$. Assume that $k+1 \in C_{m}$ and let
$x \in X_{\alpha}$ be such that $r_k(x)/q_{k+1} \geq \eps$. Then \be
\eps &\leq& \frac{x_{k+1}q_k+r_{k-1}(x)}{a_{k+1}q_k+q_{k-1}} \leq
\frac{x_{k+1}q_k+q_k+q_{k-1}}{a_{k+1}q_k} \leq
\frac{x_{k+1}+2}{a_{k+1}} \\
x_{k+1} &\geq& \eps a_{k+1}-2 = \eps_k \ee The probability of this
event is bounded away from one. For any $j\leq a_k$ we have \be \mu
\lb W_{k+1}\geq \eps_{k+1}| W_k=j] &\leq&
\frac{((1-\eps)a_{k+1}+3)\eta_k+\eta_{k+1}}{\eta_{k-1}} \\
&\leq& \frac{((1-\eps)a_{k+1}+4)\eta_k}{a_{k+1}\eta_k} \leq \delta \ee
It follows that $\mu\lb W_{k+1}\geq \eps_{k+1}|W_j=i] \geq \delta$
whenever $j\leq k$ and $i \in \{0,\ldots,a_j\}$. Given $k_0>0$, let
$k_0<k_1<k_2<\cdots <k_n$ be an incresing sequence of indices from
$C_{m}$. Then $\mu\lb W_{k_1}\geq \eps_{k_1},\ldots,W_{k_n}\geq
\eps_{k_n}] \leq \delta^n$. It follows
$$ \mu\left \{x\in X_{\alpha} \left| \frac{r_{k_i}(x)}{q_{k_{i+1}}} \geq
    \eps, 1\leq i\leq n \right \} \right. \leq \mu\lb x\in X_{\alpha}|
x_{k_i} \geq \eps_{k_i}, 1\leq i \leq n] \leq \delta^n $$ so
$\mu\{x\in X_{\alpha}: \overline{R}(x)<\frac{1}{\eps}\} = 0$ and
$\overline{R}(x) = \infty$ almost everywhere.  \eP

\begin{pro}
  If $\alpha$ have bounded coeficients in its continued fraction, then
  $\uR(x)=\rr_0$, $\oR(x)=\rr_1$ almost everywhere.
\end{pro}

\bP Proposition \ref{meze} says that $\min\uR=\uR(b)$, where
$b=(a_1,a_2,a_3,\ldots)$. We are going to prove that
$$ \mu \left\{x\in X_{\alpha} \left| \limsup_{k\to\infty}\frac{r_k(x)}{q_k}\leq
    \limsup_{k\to\infty}\frac{r_k(b)}{q_k} \right\} \right.=1.$$ Fix
$m\geq 1$. There exists an integer sequence $\{n_k\}^\infty_{k=0}$
such that $n_0>m$,\\
$n_k-n_{k-1}>m$, for $k\geq 1$ and
$$\limsup_{n\to\infty}\frac{r_n(e)}{q_n}=
\lim_{k\to\infty}\frac{r_{n_k}(e)}{q_{n_k}}.$$ For $k\in\NN$, set
$$ D_k=\{x\in X_{\alpha} \mid x_{n_k}=b_{n_k},x_{n_k-1}=b_{n_k-1},
\ldots,x_{n_k-m+1}=b_{n_k-m+1}\}.$$ and
$D=\bigcap^\infty_{j=1}\bigcup^\infty_{k=j}D_k$. We show
$\mu(D)=1$. Let $M$ be a bound for the continued fraction expansion,
so $a_k\leq M$ for every $k$. Then for any $i\leq a_k$, $j\leq
a_{k+1}$ we have
$$ \mu\lb W_{k+1}=j|W_k=i] \geq \frac{\eta_k}{\eta_k+\eta_{k+1}}
\geq \frac{1}{a_{k+1}+2} \geq \frac{1}{M+2} $$ It follows
$\mu(X_{\alpha}\setminus D_k) \leq 1 - \frac{1}{(M+2)^n}$ and
$$ \mu(D)=1-\mu \left(\bigcup^\infty_{j=1}\bigcap^\infty_{k=j}
  (X_{\alpha} \setminus D_k) \right)=1-0=1.$$ Given $x\in D$, there
exists an increasing integer sequence $\{k_j\}^\infty_{j=1}$ such that
$x\in D_{n_{k_j}}$. For each $j$, we have \be
r_{n_{k_j}}(b)-r_{n_{k_j}}(x)&=&
\sum^{n_{k_j}}_{i=0}b_{i+1}q_i-\sum^{n_{k_j}}_{i=0}x_{i+1}q_i
=  \sum^{n_{k_j}-m}_{i=0}b_{i+1}q_i-\sum^{n_{k_j}-m}_{i=0}x_{i+1}q_i \\
&=& r_{n_{k_j}-m}(b)-r_{n_{k_j}-m}(x) \leq
q_{n_{k_j}-m+1}+q_{n_{k_j}-m}-1-q_{n_{k_j}-m} \\
&\leq& q_{n_{k_j}-m+1}.  \ee Since $q_{n+2}=a_{n+2}q_{n+1}+q_n\geq
2q_n$, we get
$$ \frac{r_{n_{k_j}}(b)}{q_{n_{k_j}}}-
\frac{r_{n_{k_j}}(x)}{q_{n_{k_j}}}\leq\frac{q_{n_{k_j}-m+1}}{q_{n_{k_j}}}
\leq\frac{2^{-\left\lf\frac{m-1}{2}\right\rf}q_{n_{k_j}}}{q_{n_{k_j}}}=
2^{-\left\lf\frac{m-1}{2}\right\rf}.$$ and
$$ \limsup_{k\to\infty}\frac{r_k(x)}{q_k}\geq
\limsup_{j\to\infty}\frac{r_{k_j}(b)}{q_{k_j}}-
2^{-\left\lf\frac{m-1}{2}\right\rf} =
\limsup_{k\to\infty}\frac{r_k(b)}{q_k}-
2^{-\left\lf\frac{m-1}{2}\right\rf}.$$ It follows
$$\mu \left\{x\in X_{\alpha}
  \left|\limsup_{k\to\infty}\frac{r_k(x)}{q_k}\geq
    \limsup_{k\to\infty}\frac{r_k(b)}{q_k}-
    2^{-\left\lf\frac{m-1}{2}\right\rf} \right\} \right.=1$$ so
$\underline{R}(x) = \rr_0$ almost everywhere.  The proof for $\oR$ is
similar using the points $c$ or $d$ instead of $b$.  \eP

\begin{cor}
  Given an irrational $\alpha=\lb 0,a_1,a_2,...]$ with convergents
  $q_k$, $\gamma=\frac{\sqrt{5}+1}{2}$, set
$$ \rr_0 = \liminf_{k\to\infty} \frac{q_k}{q_{k+1}+q_k-1},\qq
\rr_1 = \limsup_{k\to\infty} \frac{q_{k+1}}{q_k} $$ Then $\rr_0\leq
\underline{R}(x) \leq \overline{R}(x) \leq \rr_1$ for every $x \in
\Sigma_{\alpha}$ and $\underline{R}(x) = \rr_0$,
$\overline{R}(x)=\rr_1$ almost everywhere.

If $\{a_k\}^\infty_{k=0}$ is unbounded, then $\rr_0=0$ and
$\rr_1=\infty$. In the case of bounded $\{a_k\}^\infty_{k=0}$,
$$\frac{1}{M+2} \leq \rr_0\leq\gamma^{-2} <\gamma\leq\rr_1\leq M+1$$
where $M=\limsup a_k$. Moreover, $\rr_1=\gamma$(
resp. $\rr_0=\gamma^{-2}$) if and only if $M=1$.
\end{cor}

\noindent{\bf Acknowledgement} I thank P.K{\accent23 u}rka from the Charles
University in Prague, whose comments were essential to the formation
of this paper. Also I'm indebted to S. Vaienti from the University
Toulon for suggesting me the problematics, during my study at Centre
de Physique Theorique in Luminy, Marseille, where I have
done the main part of the paper.\\

\end{document}